\documentclass[12pt]{article}
\usepackage{amsmath,amsfonts,calrsfs,amssymb,color}

\newtheorem{Theorem}{Theorem}[section]
\newtheorem{Definition}[Theorem]{Definition}

\newtheorem{Lemma}[Theorem]{Lemma}

\newtheorem{Remark}[Theorem]{Remark}

\makeatletter
\newif\ifmsbmloaded@

\def\loadmsbm{\msbmloaded@true
 \font\tenmsb=msbm10 scaled 1\@ptsize00
 \font\sevenmsb=msbm7 scaled 1\@ptsize00
 \font\fivemsb=msbm5 scaled 1\@ptsize00
 \alloc@8\fam\chardef\sixt@@n\msbfam
 \textfont\msbfam=\tenmsb
 \scriptfont\msbfam=\sevenmsb
 \scriptscriptfont\msbfam=\fivemsb
 }
\def\nonmatherr@#1{\errmessage%
{LateX error: \string#1\space allowed only in math mode}}
\def\Bbb{\relax\ifmmode\expandafter\Bbb@\else
 \expandafter\nonmatherr@\expandafter\Bbb\fi}
\def\Bbb@#1{{\Bbb@@{#1}}}
\def\Bbb@@#1{\fam\msbfam\relax#1}

\@addtoreset{equation}{section}

\makeatother
\loadmsbm

\def\R{\mathbb R}
\def\N{\mathbb N}

\def\E{\mathbb E}
\def\P{\mathbb P}
\def\Q{\mathbb Q}

\def\ds{\displaystyle}

\begin{document}

\title{Existence of strong solutions for stochastic porous media equation under general monotonicity conditions}
  \author{Viorel Barbu \thanks{Supported by the CEEX Project 05 of
Romanian Minister of Research.},\\
University Al. I. Cuza \\and \\
Institute of Mathematics ``Octav Mayer'', Iasi, Romania
,\\
Giuseppe Da Prato \thanks{Supported by the research program ``Equazioni di
Kolmogorov'' from the Italian
``Ministero della Ricerca Scientifica e Tecnologica''},\\
 Scuola Normale Superiore
di Pisa, Italy\\
 and \\
Michael R\"ockner \thanks{Supported  by the SFB-701 and the
BIBOS-Research Center.} \\
Faculty of Mathematics, University of Bielefeld, Germany\\
and\\
Department of Mathematics and Statistics, Purdue University,\\  U. S. A.}

\maketitle

{\bf Abstract}. One proves existence and uniqueness of strong solutions to
  stochastic porous media equations under minimal monotonicity conditions
on the nonlinearity. In particular, we do not assume continuity of the drift or any growth
condition at infinity.



\section{Introduction}

This work is concerned with existence and uniqueness of solutions to
  stochastic porous media equations
\begin{equation}
\label{e1.1}
\left\{\begin{array}{l}
dX(t)-\Delta\Psi(X(t))dt=B(X(t))dW(t)\quad\mbox{\rm in}\;(0,T)\times \mathcal O:=Q_T,\\
\Psi(X(t))=0\quad\mbox{\rm on}\;(0,T)\times \partial\mathcal O:=\Sigma_T,\\
X(0)=x\quad\mbox{\rm in}\; \mathcal O,
\end{array}\right.
\end{equation}
where $\mathcal O$ is an open, bounded domain of $\R^d$, $d\ge 1,$ with smooth boundary $\partial\mathcal O$,
$W(t)$ is a cylindrical Wiener process on $L^2(\mathcal O)$, 
  $B:H\to L(L^2(\mathcal O),L^2(\mathcal O))$ is a Lipschitz continuous operator to be precised below
and  $H:=H^{-1}(\mathcal O)$.
The function $\Psi:\R\to \R$ (or more generally the multivalued function $\Psi:\R\to 2^\R$) is a maximal monotone
graph in $\R\times \R$.\bigskip

Existence results for equation \eqref{e1.1} were obtained in \cite{9} (see also \cite{3},\cite{4}) in the  special 
case $B=\sqrt{Q}$, with $Q$ linear nonnegative, Tr $Q<+\infty$ and $\Psi\in C^1(\R)$ satifying the growth condition
\begin{equation}
\label{e1.2}
k_3+k_1|s|^{r-1}\le \Psi'(s)\le k_2(1+|s|^{r-1}), \quad s\in \R,
\end{equation}
where $k_1,k_2>0$, $k_3\in \R$, $r>1.$

Under these growth conditions on   $\Psi$, equation \eqref{e1.1} covers many important models of
dynamics of an ideal gas in a porous medium (see e.g. \cite{1}) but excludes, however, other significant physical
models such as plasma fast diffusion (\cite{5}) which arises for $\Psi(s)=\sqrt{s}$, phase transitions or dynamics
of saturated underground water flows (Richard's equation). In the later case  multivalued monotone
graphs $\Psi$ might appear (see \cite{12}). Recently in \cite{15} (see also \cite{14}) 
the existence results of \cite{9} were
extended to the case of monotone nonlinearities $\Psi$ such that $s\mapsto s\Psi(s)$ is (comparable to) a $\Delta_2$-regular Young function
(cf. assumption (A1) in \cite{15}) thus including the fast diffusion model.
As a matter of fact, in the line of the classical work of N. Krylov and B. Rozovskii \cite{10} the approach used in
\cite{15} is a variational one i.e. one considers the stochastic equation \eqref{e1.1} in a duality setting induced
by a functional triplet $V\subset H\subset V'$ and this requires to find appropriate spaces $V$ and $H$. This was done in 
\cite{15} in an elaborate way even with $\Delta$ replaced by very general (not necessarily differential) operators
$L$.

The method we use here is quite different and essentially  an $L^1$-approach relying on   weak compacteness
techniques in $L^1(Q_T)$ via the Dunford-Pettis theorem which involve minimal growth assumptions on $\Psi$.
Restricted to single valued continuous functions $\Psi$ the main result, Theorem \ref{t2.2} below, gives existence
and uniqueness of solutions only assuming that $\lim_{s\to+\infty}\Psi(s)=+\infty$,
$\lim_{s\to-\infty}\Psi(s)=-\infty$, $\Psi$ monotonically increasing and
\begin{equation}
\label{e1.3}
\limsup_{|s|\to +\infty}\;\frac{\int_0^{-s}\Psi(t)dt}{\int_0^{s}\Psi(t)dt}<+\infty.
\end{equation}
We note that the assumptions on $\Psi$ in \cite{15}) imply our assumptions.
 In this sense, under assumption (H2) below in  the noise, the results on this paper extend those in \cite{15} in case
$L=\Delta$    if $\mathcal O$ is bounded and if the coefficients do not depend on $(t,\omega)$. 
The latter two were not assumed in \cite{15}. On the other hand a growth condition on $\Psi$ is imposed in
\cite{15} (cf. \cite[Lemma 3.2]{15}) which is not done here. Another main progress of this paper  is that
$\Psi$ is no longer assumed to be continuous, it might be multivalued and with exponential growth to $\pm \infty$ (for instance of the form
$\exp{(a|x|^p)}$).
We note  that \eqref{e1.3} is not a growth condition at $+\infty$ but a kind of symmetry
condition about the  behaviour
 of $\Psi$ at $\pm \infty$. If $\Psi$ is a maximal monotone graph with   potential $j$ (i.e.
$\Psi=\partial j$) then \eqref{e1.3} takes the form (see Hypothesis $(H3)$ below)
$$
\limsup_{|s|\to +\infty}\;\frac{j(-s)}{j(s)}<+\infty.
$$
Anyway this condition is automatically satified for even monotonically increasing functions $\Psi$ or e.g. if a
condition of the form \eqref{e1.2} is satisfied. We note, however, that because of our very general conditions on $\Psi$ the solution of
\eqref{e1.1} will be pathwise only weakly continuous in $H$. It seems impossible to prove strong continuity.

\subsection{Notations}

$\mathcal O$ is a bounded open subset of $\R^d$, $d\ge 1$ with smooth boundary $\partial\mathcal O$. We set
$$
Q_T=(0,T)\times \mathcal O,\quad\Sigma_T=(0,T)\times \partial\mathcal O.
$$
Moreover $L^p(\mathcal O)$, $L^p(Q_T)$, $p\ge 1$, are standard $L^p$- function spaces and $H^1_0(\mathcal O)$,
$H^k(\mathcal O)$ are Sobolev spaces on $\mathcal O$. Denote by $H:=H^{-1}(\mathcal O)$ the dual of $H^1_0(\mathcal
O)$ with the norm and the scalar product
$$
|u|_{-1}:=(A^{-1}u,u)^{1/2},\quad \langle  u,v \rangle_{-1}=(A^{-1}u,v),
$$
where $(\cdot,\cdot)$ is the pairing between $H^1_0(\mathcal O)$ and $H^{-1}_0(\mathcal O)$ which coincides with the
scalar product of $L^2(\mathcal O)$. We have denoted by $A$ the Laplace operator with Dirichlet homogeneous boundary
conditions, i.e.
\begin{equation}
\label{e1.4}
Au=-\Delta u,\quad u\in D(A)=H^{2}(\mathcal O)\cap H^{1}_0(\mathcal O).
\end{equation}
Given a Hilbert space $U$, the norm of $U$ will be denoted by $|\cdot |_U$ and the scalar product by
$(\cdot,\cdot)_U$. By $C([0,T];U)$ we shall denote the space of $U$-valued continuous functions on $[0,T]$ and by
$C^w([0,T];U)$  the space of weakly continuous functions from $[0,T]$ to $U$.

Given two Hilbert spaces $U$ and $V$ we shall denote by $L(U,V)$ the space of linear continuous operators from $U$
to $V$ and by $L_{HS}(U,V)$ the space of Hilbert-Schmidt operators $F:U\to V$ with the norm
\begin{equation}
\label{e1.5}
\|F\|_{L_{HS}(U,V)}:=\left(\sum_{i=1}^\infty|Fe_i|^2_V\right)^{1/2},
\end{equation}
where $\{e_i\}$ is an orthonormal basis in $U$.

If $j:\R\to (-\infty,+\infty]$ is a lower semicontinuous convex function we denote by $\partial j:\R\to 2^{\R}$
the subdifferential of $j$, i.e.
$$
\partial j(y)=\{\theta\in \R:\;j(y)\le j(z)+\theta(y-z),\;\forall\;z\in \R\}
$$
and by $j^*$ the conjugate of $j$ (the Legendre transform of $j$),
$$
j^*(p)=\sup\{py-j(y):\;y\in \R\}.
$$
We recall that $\partial j^*=(\partial j)^{-1}$  (see e.g. \cite{2}, \cite{6}),
\begin{equation}
\label{e1.6}
j(y)+j^*(p)=py\quad\mbox{\rm if and only if}\;p\in \partial j(y)
\end{equation}
and
\begin{equation}
\label{e1.7}
j(u)+j^*(p)\ge pu\quad\mbox{\rm for all}\;p,u\in \R.
\end{equation}
Moreover  $\Psi:=\partial j$ is maximal monotone, i.e.
$$
(y_1-y_2)(p_1-p_2)\ge 0\quad\mbox{\rm for all}\;p_i\in \partial j(y_i),\; i=1,2
$$
and $R(1+\partial j)=\R$.

Given a multivalued function $\Phi:\R\to 2^\R$ we shall denote by $D(\Phi)=\{u\in \R:\;\Phi(u)\neq \varnothing\}$
the domain of $\Phi$ and by $R(\Phi)=\{v:\;v\in \Phi(u),\;u\in D(\Phi)\}$ its range.

Given a maximal monotone graph $\Psi:\R\to 2^\R$ there is a unique lower semicontinuous convex function
$j:\R\to (-\infty,+\infty]$ such that $\Psi:=\partial j$. The function $j$ is unique up to an additive constant and   called the {\em potential} of
$\Psi$.

For any maximal monotone graph $\Psi$ and $\lambda>0$   by
$$
\Psi_\lambda=\frac1\lambda\;(1-(1+\lambda\Psi)^{-1})\in \Psi(1+\lambda\Psi)^{-1}
$$
we denote the {\em Yosida} approximation of $\Psi$.  Here $1$ stands for the identity function.
$\Psi_\lambda$ is   Lipschitzian and
monotonically increasing.

\section{The main result}

\subsection{Hypotheses}

\noindent $(H_1)$  $W(t)${\it  is a cylindrical Wiener process on $L^2(\mathcal O)$ defined by
\begin{equation}
\label{e2.1}
W(t)=\sum_{k=1}^\infty\beta_k(t)e_k,
\end{equation}
where $\{\beta_k\}$ is a sequence of mutually independent Brownian motions on a filtered probability space
$(\Omega,\mathcal F,\{\mathcal F_t\}_{t\ge 0},\P)$, with right continuous filtration and  $\{e_k\}$ is an orthonormal basis in
$L^2(\mathcal O)$. To be more specific $\{e_k\}$ will be chosen as the normalized sequence of eigenfunctions of 
the operator $A$, hence $e_k\in L^p(\mathcal O)$ for all $k\in \N,p\ge 1$}.\bigskip

\noindent $(H_2)$ {\it  $B$ is Lipschitzian from $H=H^{-1}(\mathcal O)$ to $L_{HS}(L^2(\mathcal O),D(A^\gamma))$
where $\gamma>d/2$.}\bigskip

\noindent $(H_3)$ {\it $\Psi:\R\to 2^\R$ is a maximal monotone graph on $\R\times \R$ such that $0\in \Psi(0),$
\begin{equation}
\label{e2.2}
D(\Psi)=\R,\quad R(\Psi)=\R
\end{equation}
and
\begin{equation}
\label{e2.3}
\limsup_{|s|\to +\infty}\;\frac{j(-s)}{j(s)}<+\infty.
\end{equation}}
Here $j:\R\to \R$ is the potential of $\Psi$ i.e. $\partial j=\Psi$, which under assumption \eqref{e2.2}
is a continuous convex function. Since $0\in \Psi(0)$, by definition we have $j(0)=\inf j$. Hence subtracting $j(0)$
we can take $j$ such that $j(0)=0$ and $j\ge 0$, hence $j^*\ge j^*(0)=0$.
It should be recalled (see e.g. \cite{2},\cite{6}) that the condition $R(\Psi)=\R$ is equivalent to
\begin{equation}
\label{e2.4}
j(y)<\infty\;\;\forall\;y\in \R,\quad \lim_{|y|\to\infty}\frac{j(y)}{|y|}=+\infty
\end{equation}  
and that the condition $D(\Psi)=\R$ is equivalent to
\begin{equation}
\label{e2.4'}
j^*(y)<\infty\;\;\forall\;y\in \R,\quad \lim_{|y|\to\infty}\frac{j^*(y)}{|y|}=+\infty
\end{equation}   
Hypothesis $(H_3)$ automatically holds if $\Psi$ is a monotonically increasing, continuous function on $\R$ satisfying
condition \eqref{e1.3} and
$$
\lim_{s\to +\infty}\;\Psi(s)=+\infty,\quad \lim_{s\to -\infty}\;\Psi(s)=-\infty.
$$
In particular, it is satisfied by functions $\Psi$ satisfying
  \eqref{e1.2} for $r>0$ or more generally by those satisfying
assumption $(A1)$   in \cite{15}.\bigskip

We need some more notations. Given a Banach space $Z$ we shall denote by
$$
C_W([0,T];Z)=C([0,T];L^2(\Omega,\mathcal F,\P;Z))
$$
the space of all continuous adapted stochastic processes which are  mean square continuous.
The space
$$
L^2_W([0,T];Z)=L^2([0,T];L^2(\Omega,\mathcal F,\P:Z))
$$
is similarly defined (see e.g. \cite{7}, \cite{8}).\bigskip

\begin{Definition}
\label{d2.1}
 An adapted process $X \in C_W([0,T];H)\cap L^1((0,T)\times\mathcal O\times \Omega)$, such that $X\in C^w([0,T],H), \P$-a.s., is said to be a
strong solution to equation \eqref{e1.1} if there exists a  
  process $\eta\in L^1((0,T)\times\mathcal O\times \Omega)$ such that
\begin{equation}
\label{e2.5}
\eta(t,\xi)\in \Psi(X(t,\xi)),\quad\mbox{\it a.e.}\;(t,\xi)\in Q_T,\;\P\mbox{\it -a.s.}
\end{equation}
\begin{equation}
\label{e2.6}
\int_0^{\bullet}\eta(s)ds\in C^w([0,T];H^1_0(\mathcal O)),\ 
\end{equation}
\begin{equation}
\label{e2.7}
X(t)-\Delta\int_0^t\eta(s)ds=x+\int_0^tB(X(s))dW(s),\quad\forall\;t\in [0,T],\;\P\mbox{\it -a.s.}
\end{equation}
\begin{equation}
\label{e2.8}
j(X),\;j^*(\eta)\in L^1((0,T)\times\mathcal O\times \Omega).
\end{equation}
(Here $\int_0^{t}\eta(s)ds$ is initially defined as on $L^1(\mathcal O)$-valued Bochner integral).
Of course, if $\Psi$ is single valued \eqref{e2.5}-\eqref{e2.7} reduce to
\begin{equation}
\label{e2.9}
\int_0^\bullet\Psi(X(s))ds\in C^w([0,T];H^1_0(\mathcal O)), 
\end{equation}
and
\begin{equation} 
X(t)-\Delta\int_0^t\Psi(X(s))ds=x+\int_0^tB(X(s))dW(s),\quad\forall\;t\in [0,T],\;\P\mbox{\it -a.s.}
\end{equation}
\end{Definition}
We note that $X$ is as in Definition \ref{d2.1} is automatically predictable.

Theorem \ref{t2.2} below is the main result of this work.
\begin{Theorem}
\label{t2.2}
Under Hypotheses $(H_1)$, $(H_2)$, $(H_3)$, for each $x\in H$ there is a unique strong solution $X=X(t,x)$ to
equation \eqref{e1.1}. Moreover, the following estimate holds
\begin{equation}
\label{e2.11}
\E|X(t,x)-X(t,y)|^2_{-1}\le C|x-y|^2_{-1},\quad\mbox{\it for all}\;t\ge 0,
 \end{equation}
where $C$ is independent of $x,y\in H$.
\end{Theorem}
Theorem \ref{t2.2} will be proved in section 4 via fixed point arguments. Previously, we shall establish
in section 3 the
existence  of solutions for the equation
\begin{equation}
\label{e2.12}
\left\{\begin{array}{l}
dY(t)-\Delta\Psi(Y(t))dt=G(t)dW(t)\quad\mbox{\rm in}\; Q_T,\\
\Psi(Y(t))=0\quad\mbox{\rm on}\;\Sigma_T,\\
Y(0)=x\quad\mbox{\rm in}\; \mathcal O,
\end{array}\right.
\end{equation}
where $G:[0,T]\to L_{HS}(L^2(\mathcal O),D(A^\gamma))$ is a predictable process such that
\begin{equation}
\label{e2.13}
\E\int_0^T\|G(t)\|^2_{L_{HS}(L^2(\mathcal O),D(A^\gamma))}dt<+\infty
\end{equation}
and $\gamma>d/2$. By $GdW$ we mean of course
$$
GdW=\sum_{k=1}^\infty Ge_k d\beta_k.
$$
By a solution of \eqref{e2.12} we shall mean   an adapted process $Y$ satisfying along with
$\eta\in L^1((0,T)\times\mathcal O\times \Omega)$  conditions \eqref{e2.5}-\eqref{e2.8} where $B(X)$ is replaced by
$G$.
\begin{Theorem}
\label{t2.3}
Under Hypotheses $(H_1)$,  $(H_3)$, \eqref{e2.13}, for each $x\in H$ there is a unique strong solution $Y=Y_G(t,x)$
to equation \eqref{e2.12} in the sense of Definition $\ref{d2.1}$. Moreover, the following estimate holds
\begin{equation}
\label{e2.14}
\begin{array}{l}
\E|Y_{G_1}(t,x)-Y_{G_2}(t,y)|^2_{-1}\le |x-y|^2_{-1}\\
\\
\ds+\E\int_0^t\|G_1(s)-G_2(s)\|^2_{L_{HS}(L^2(\mathcal O),H)}ds
 ,\quad\mbox{\it for all}\;t\ge 0,
\end{array}
 \end{equation}
for all $x,y\in H$ and $G_1,G_2$ satisfying \eqref{e2.13}.
\end{Theorem}
\begin{Remark}
\label{r2.4}
\em It should be noted that assumption (H2) excludes the case of equation \eqref{e1.1} with covariance operator $B$ of the
form $B(X)=X$ i.e. the case of multiplicative noise. However such an equation can be approximated
taking $B(X)=X*\rho_\epsilon$ ($\rho$ is a mollifier in $\mathcal O$) or $B(X)=(1+\epsilon A)^{-\delta} X,\;
\epsilon>0.$
\end{Remark}
\begin{Remark}
\em Assumption  $(H_3)$ for example allows monotonically increasing functions $\Psi$ which are continuous from the right on $\R$
and have a finite number of jumps  $r_1,r_2,... ,r_N$. However in this case one must
fill the jumps by replacing the function $\Psi$ by the maximal
monotone (multivalued ) graph $\tilde\Psi(r) = \Psi(r)$ for $r$ different from
$r_i$ and $\tilde\Psi(r_i) = [\Psi(r_i)-\Psi(r_{i-1}-0)]$. Such a situation might arise in modelling of underground water
flows (see e.g. \cite{12}). In this case $\Psi$ is the diffusivity function and \eqref{e1.1} reduces to Richard's
equation.
It must be also said   that Theorems \ref{t2.2} and \ref{t2.3} have
natural extensions to equations of the form
\begin{equation}
\label{e2.15}
dX(t)-\Delta\Psi(X(t))dt+\Phi(X(t))dt=B(X(t))dW(t),
\end{equation}
where $\Phi$ is  a suitable monotonically increasing and continuous function (see \cite{15}).
 Also as in \cite{15} one might consider the case where
$\Psi=\Psi(X,\omega),\;\omega\in \Omega,$ but we do not go into details, here.
We  also note that assumption $D(\Psi)=\R$ in Hypothesis $(H_3)$  excludes a situation of the following type
\begin{equation}
\label{e2.16}
\Psi(s)=\left\{\begin{array}{l}
\Psi_1(s)\quad\mbox{\rm for }\;s<s_0,\; \Psi(s_0)=(0,+\infty),\\
=\varnothing\quad\mbox{\rm for }\;s>s_0,
\end{array}\right.
\end{equation}
where  $\Psi$ is a continuous monotonically increasing function such that $\Psi_1(-\infty)$ $=-\infty.$
In this case problem
\eqref{e1.1} reduces to a stochastic variational inequality and it is relevant in the description of saturation processes in
infiltration.
An analysis similar to that to be developped
below shows that in this case in Definition \ref{d2.1} the solution is no more an $L^1$- function but a bounded measure on $Q_T$.
We expect to give details in a later paper.

Another situation of interest covered by our assumptions (see also \cite{15}) is that of logarithmic diffusion equations arising in plasma
physics see e.g. \cite{13}. In this case $\Psi(s)=\log(\mu+|s|)$ sign$(s)$.
\end{Remark}

\section{Proof of Theorem \ref{t2.3}}

For every $\lambda>0$ consider the approximating equation
\begin{equation}
\label{e3.1}
\left\{\begin{array}{l}
dX_\lambda(t)-\Delta(\Psi_\lambda(X_\lambda(t))+\lambda X_\lambda(t))dt=G(t)dW(t)\quad\mbox{\rm in}\;(0,T)\times
\mathcal O:=Q_T,\\
\Psi_\lambda(X_\lambda(t))+\lambda X_\lambda(t))=0\quad\mbox{\rm on}\;(0,T)\times \partial\mathcal O,\\
X_\lambda(0)=x\quad\mbox{\rm in}\; \mathcal O,
\end{array}\right.
\end{equation}
which has a unique solution $X_\lambda\in C_W([0,T];H)$ such that 
$$X_\lambda,\Psi_\lambda(X_\lambda)\in
L^2_W(0,T;H^1_0(\mathcal O)).$$
Indeed, setting $y_\lambda(t)=X_\lambda(t)-W_G(t)$ where $W_G(t)=\int_0^tG(s)dW(s)$, we may rewrite \eqref{e3.1} as a
  random equation
\begin{equation}
\label{e3.2}
\left\{\begin{array}{l}
y'_\lambda(t)-\Delta\widetilde{\Psi}_\lambda(y_\lambda(t)+W_G(t))=0\quad\P\mbox{\rm -a.s. in }\; Q_T,\\
\widetilde{\Psi}_\lambda(y_\lambda(t)+W_G(t)) =0 \quad\mbox{\rm on}\;(0,T)\times
\partial\mathcal O,\\ 
y_\lambda(0)=x\quad\mbox{\rm in}\; \mathcal O,
\end{array}\right.
\end{equation}
where $\widetilde{\Psi}_\lambda(y)= \Psi_\lambda(y)+\lambda y,\;\lambda>0.$ Note that
$\widetilde{\Psi}_\lambda(0)=0$.

For each $\omega\in \Omega$ the operator $\Gamma(t):H^1_0(\mathcal O)\to H^{-1}(\mathcal O)$ defined by
$$
\Gamma(t)y=-\Delta\widetilde{\Psi}_\lambda(y+W_G(t)),\quad y\in H^1_0(\mathcal O),
$$
is continuous, monotone and  coercive, i.e.
$$(\Gamma(t)y,y)\ge \lambda |y+W_G(t)|^2_{H^1_0(\mathcal O)}-(\Gamma(t)y,W_G(t))\ge
\frac\lambda2\;|y|^2_{H^1_0(\mathcal O)}- C_\lambda|W_G(t)|^2_{H^1_0(\mathcal O)}.
$$
Then by classical existence theory for nonlinear equations (see e.g. \cite{11}) equation \eqref{e3.2}
has a unique solution  
$$y_\lambda\in C([0,T];L^2(\mathcal O))\cap L^2(0,T;H^{1}_0(\mathcal O))$$
with $y'_\lambda\in L^2(0,T;H^{-1}(\mathcal O))$. The function $X_\lambda(t)=y_\lambda(t)+W_G(t)$
is of course an adapted process because the solution $y_\lambda$ to equation \eqref{e3.2} is a continuous function
of $W_G$ and so it satifies the requested condition.

\subsection{A-priori estimates}
From now on we shall fix $\omega\in \Omega$ and work with the corresponding solution $y_\lambda$ to
\eqref{e3.2}. We have
\begin{equation}
\label{e3.3}
\begin{array}{l}
\ds
\frac12\;\frac{d}{dt}\;|y_\lambda(t)|^2_{-1}+(\widetilde{\Psi}_\lambda(y_\lambda(t)+W_G(t)),y_\lambda(t)+W_G(t))\\
\\
=(\widetilde{\Psi}_\lambda(y_\lambda(t)+W_G(t)),W_G(t)),
\end{array} 
\end{equation}
which is equivalent to
\begin{equation}
\label{e3.4}
\begin{array}{l}
\ds
\frac12\;\frac{d}{dt}\;|y_\lambda(t)|^2_{-1}+(\Psi_\lambda(y_\lambda(t)+W_G(t)),y_\lambda(t)+W_G(t))\\
\\
=-\lambda (y_\lambda(t),y_\lambda(t)+W_G(t))+( \Psi_\lambda(y_\lambda(t)+W_G(t)),W_G(t)).
\end{array} 
\end{equation}
Now  set $j_\lambda(u)=\int_0^u\Psi_\lambda(r)dr$ and denote by $j^*_\lambda$  the conjugate of $j_\lambda$. 
 By \eqref{e1.6} we have
$$
j^*_\lambda (\Psi_\lambda(y_\lambda(t)+W_G(t)))
 +j_\lambda (y_\lambda(t)+W_G(t))=\Psi_\lambda(y_\lambda(t)+W_G(t))(y_\lambda(t)+W_G(t)).
$$
Substituting this identity into \eqref{e3.4}
yields
\begin{equation}
\label{e3.5}
\begin{array}{l}
\ds
\frac12\;|y_\lambda(t)|^2_{-1}+ \int_0^t\int_\mathcal
O(j_\lambda(y_\lambda(s)+W_G(s))+j^*_\lambda(\Psi_\lambda(y_\lambda(s)+W_G(s)))d\xi ds\\
\\\
\ds=
\frac12\;|x|^2_{-1}+ \int_0^t\int_\mathcal
O(\Psi_\lambda(y_\lambda(s)+W_G(s)) W_G(s))d\xi ds\\
\\
\ds-\lambda\int_0^t\int_\mathcal
Oy_\lambda(s) (y_\lambda(s)+W_G(s))d\xi ds,
\end{array} 
\end{equation}
Since $j_\lambda$ is the Moreau approximation of $j$
$$
j_\lambda(u)=\min\left\{j(v)+\frac1{2\lambda}\;|u-v|^2:\; v\in H\right\},
$$
we have (recall that the minimum is attained at $v=(1+\lambda\Psi)^{-1}u$)
\begin{equation}
\label{e3.6}
j_\lambda(u)=j((1+\lambda\Psi)^{-1}u)+\frac1{2\lambda}\;|u-(1+\lambda\Psi)^{-1}u|^2,\quad u\in \R.
\end{equation}
We  now set
\begin{equation}
\label{e3.7}
z_\lambda=(1+\lambda\Psi)^{-1}(y_\lambda+W_G),\quad\eta_\lambda=\Psi_\lambda(y_\lambda+W_G).
\end{equation}
Then, using \eqref{e3.6} and the fact that $j^*_\lambda\ge j^*$ for all $\lambda>0$, we see by \eqref{e3.5} that
\begin{equation}
\label{e3.8}
\begin{array}{l}
\ds
\frac12\;|y_\lambda(t)|^2_{-1}+ \int_0^t\int_\mathcal
O(j(z_\lambda(s))+j^*(\eta_\lambda(s))d\xi ds\\
\\
\ds +\frac1{2\lambda}\;\int_0^t\int_\mathcal O
(y_\lambda(s)+W_G(s)-z_\lambda(s))^2d\xi ds\\
\\
\ds\le
\frac12\;|x|^2_{-1}+ \int_0^t\int_\mathcal
O\eta_\lambda(s)W_G(s)d\xi ds  -\lambda\int_0^t\int_\mathcal
O  y_\lambda(s)(y_\lambda(s)+W_G(s))d\xi ds.
\end{array} 
\end{equation}
We now estimate the first integral from the right hand side of \eqref{e3.8} as follows
\begin{equation}
\label{e3.9}
 \left| \int_0^t\int_\mathcal
O  \eta_\lambda(s)W_G(s) d\xi ds\right|\le  \delta \int_0^t\int_\mathcal
O  |\eta_\lambda(s)|   d\xi ds,
\end{equation}
where $\delta:=\sup_{s\in [0,T]}|W_G(s)|_{L^\infty(\mathcal O)}<+\infty$.
We note that by assumption \eqref{e2.13} and since $\gamma>d/2$ it follows 
by Sobolev embedding that $W_G(\cdot)$ has  continuous sample paths in    $D(A^\gamma)\subset
L^\infty(\mathcal O)$ and so   $\delta$  is indeed finite.

Substituting \eqref{e3.9} in \eqref{e3.8} yields
$$
\begin{array}{l}
\ds
\frac12\;|y_\lambda(t)|^2_{-1}+ \int_0^t\int_\mathcal
O(j(z_\lambda(s))+j^*(\eta_\lambda(s)d\xi ds\\
\\
\ds +\frac1{2\lambda}\;\int_0^t\int_\mathcal O
(y_\lambda(s)+W_G-z_\lambda(s))^2 d\xi ds\\
\\
\ds\le
\frac12\;|x|^2_{-1}+\delta \int_0^t\int_\mathcal
O  |\eta_\lambda(s)|   d\xi ds  -\lambda\int_0^t\int_\mathcal
O  y_\lambda(s)(y_\lambda(s)+W_G(s)) d\xi ds.
\end{array} 
$$
Since  
$$
-y_\lambda(s)(y_\lambda(s)+W_G(s))\le -\frac12\;|y_\lambda(s)|^2+\frac12\;W_G^2(s),
$$
we find 
\begin{equation}
\label{e3.10}
\begin{array}{l}
\ds
\frac12\; |y_\lambda(t)|^2_{-1}+ \int_0^t\int_\mathcal
O(j(z_\lambda(s))+j^*(\eta_\lambda(s))d\xi ds+\frac{\lambda}{2}\;\int_0^t\int_\mathcal O |y_\lambda(s)|^2 d\xi ds
\\
\\
\ds+\frac1{2\lambda}\;\int_0^t\int_\mathcal O (y_\lambda(s)+W_G(s)-z_\lambda(s))^2 d\xi ds
\\
\\
\ds\le \left(\frac12\;|x|^2_{-1}+ \delta\int_0^t\int_\mathcal
O|\eta_\lambda(s)| d\xi ds+ \frac\lambda2\int_0^t\int_\mathcal
O  W_G^2(s) d\xi ds\right),\quad t\in [0,T].
\end{array} 
\end{equation}
On the other hand, we recall that condition $D(\Psi)=\R$ is equivalent with
\begin{equation}
\label{e3.11}
j^*<\infty\;\;\mbox{\rm and}\;\;\lim_{|p|\to\infty}\frac{j^*(p)}{|p|}=+ \infty.
\end{equation}
So, there exists $N=N(\omega)$ such that
$$
|\eta_\lambda(s)|> N \Rightarrow j^*(\eta_\lambda(s))>2C\delta|\eta_\lambda(s)|.
$$
Consequently, for $C>|Q_T|$  we have
$$
\begin{array}{lll}
\ds   \int_0^t\int_\mathcal
O|\eta_\lambda(s)| d\xi ds &=&\ds\int\int_{|\eta_\lambda(s)|> N}|\eta_\lambda(s)| d\xi ds+
\int\int_{|\eta_\lambda(s)|\le N}|\eta_\lambda(s)| d\xi ds
\\
\\
  &\le &\ds\frac1{2C\delta}\; \int_0^t\int_\mathcal Oj^*(\eta_\lambda(s)) d\xi ds+NC\delta.
\end{array}
$$
Substituting this into \eqref{e3.10}, since $j\ge 0$, we obtain the estimate
\begin{equation}
\label{e3.12}
\begin{array}{l}
\ds
 \frac12\;|y_\lambda(t)|^2_{-1}+ \int_0^t\int_\mathcal
O(j(z_\lambda(s))+j^*(\eta_\lambda(s)) d\xi ds\\
\\\
\ds+\frac1{2\lambda}\;\int_0^t\int_\mathcal O (y_\lambda+W_G-z_\lambda)^2 d\xi ds \le C_1 (1+|x|^2_{-1}) ,\quad t\in
[0,T],
\end{array} 
\end{equation}
which implies
\begin{equation}
\label{e3.13}
\int_0^t\int_\mathcal
O(j(z_\lambda(s))+j^*(\eta_\lambda(s)) d\xi ds\le C_1 (1+|x|^2_{-1}),
\end{equation}
and
\begin{equation}
\label{e3.14}
\int_0^t\int_\mathcal O (y_\lambda+W_G-z_\lambda)^2d\xi ds \le 2\lambda C_1 (1+|x|^2_{-1}),
\end{equation}
where $C_1$ is a  suitable random constants.

\subsection{Convergence for $\lambda\to 0$}
Since (by \eqref{e2.4} and \eqref{e2.4'})
\begin{equation}
\label{e3.14'}
\lim_{|u|\to \infty} j(u)/|u| =\infty,\quad \lim_{|u|\to \infty} j^*(u)/|u| =\infty,
\end{equation} 
we deduce from \eqref{e3.13} that the sequences $\{z_\lambda\}$ and $\{\eta_\lambda\}$
are bounded and equi-integrable in $L^1(Q_T)$. 
Then by the Dunford-Pettis theorem the sequences $\{z_\lambda\}$ and $\{\eta_\lambda\}$  are weakly compact
in $L^1(Q_T)$. Hence on a subsequence, again denoted by $\lambda$, we have
\begin{equation}
\label{e3.15}
z_\lambda\to z,\;\;\eta_\lambda\to \eta\quad\mbox{\rm weakly in }\;L^1(Q_T)\;\mbox{\rm as}\;\lambda\to 0.
\end{equation}
Moreover, by \eqref{e3.12}, \eqref{e3.14} we see that $z=y+W_G$ where
\begin{equation}
\label{e3.16}
y_\lambda\to y \quad\mbox{\rm weakly$^*$ in }\;L^\infty(0,T;H)\;\mbox{\rm and weakly in }\;L^1(Q_T).
\end{equation}
Note also that by \eqref{e3.2} we have for every $t\in [0,T]$
\begin{equation}
\label{e3.17}
y_\lambda(t)-\Delta\left(\int_0^t(\eta_\lambda(s)+
\lambda (y_\lambda(s)+W_G(s)))ds\right)=x
\end{equation}
and so the sequence $\left\{\int_0^\bullet(\eta_\lambda(s)+\lambda y_\lambda(s))ds\right\}$
is bounded in $L^\infty(0,T;H^1_0(\mathcal O))$. Hence, selecting a further subsequence if necessary (see
\eqref{e3.10}),
we have
\begin{equation}
\label{e3.18}
\lim_{\lambda\to 0}\int_0^\bullet(\eta_\lambda(s)+\lambda y_\lambda(s))ds=\int_0^\bullet\eta (s) 
 ds\;\;\mbox{\rm weakly$^*$ 
in }\;L^\infty(0,T;H^1_0(\mathcal O)).
\end{equation}
So, by \eqref{e3.17} we find
\begin{equation}
\label{e3.19}
y(t)+A\int_0^t\eta(s)ds=x\quad\;\mbox{\rm a.e.}\;t\in [0,T].
\end{equation}
Since
$$
\int_0^\bullet\eta (s) ds\in C([0,T];L^1(\mathcal O))\cap L^\infty(0,T;H_0^1(\mathcal O)),
$$
$t\mapsto \int_0^t\eta(s)ds$ is weakly continuous in $H_0^1(\mathcal O)$, hence so is 
$t\mapsto A\int_0^t\eta(s)ds$ in $H$. So, defining
\begin{equation}
\label{e3.20}
\tilde{y}(t):=-A\int_0^t\eta(s)ds+x,\quad t\in [0,T],
\end{equation}
$\tilde{y}$ is an $H$-valued weakly continuous version of $y$. Furthermore, we claim that  for
$\lambda\to 0$
$$
y_\lambda(t)\to \tilde{y}(t) \quad\mbox{\rm weakly  in }\;H,\quad\forall\;t\in [0,T].
$$
Indeed, since $\eta_\lambda\to \eta$ weakly    in $L^1(Q_T)$ and $\lambda(y_\lambda+W_G)\to 0$
weakly  in    $L^1(Q_T)$ (since it even converges strongly in $L^2(Q_T)$ to zero by \eqref{e3.10}), it follows that
for every $t\in [0,T]$
$$
\int_0^t(\eta_\lambda(s)+\lambda(y_\lambda(s)+W_G(s)))ds\to \int_0^t\eta(s)ds
 \quad\mbox{\rm weakly  in }\;L^1(\mathcal O).
$$
Hence by \eqref{e3.17}) and the definition of $\tilde{\eta}$ we obtain that for every $t\in [0,T]$
$$
(-\Delta)^{-1}y_\lambda(t)\to (-\Delta)^{-1}\tilde{y}(t) \quad\mbox{\rm weakly  in }\;L^1(\mathcal O).
$$
Since $y_\lambda(t),\; \lambda> 0,$ are bounded in $H$ by \eqref{e3.12}, the above immediately implies the claim.

From now on we always consider this particular version $\tilde{y}$ of $y$ defined in \eqref{e3.20}).
For simplicity we denote it again by $y$; so we have
$$
y_\lambda(t)\to  y(t) \quad\mbox{\rm weakly  in }\;H,\quad\forall\;t\in [0,T].
$$
We can also rewrite  
equation \eqref{e3.20}  as
\begin{equation}
\label{e3.21}
y_t(t)-\Delta \eta (t) =0\quad\mbox{\rm in }\;\mathcal D\,'(Q_T),\quad y(0)=x.
\end{equation}
Now we are going to show that
\begin{equation}
\label{e3.22}
\eta(t,\xi)\in\Psi(y(t,\xi)+W_G(t,\xi))\quad\mbox{\rm a.e. }\;(t,\xi)\in Q_T.
\end{equation}
For this we shall need 
the following inequality 
\begin{equation}
\label{e3.23}
\liminf_{\lambda\to 0} \int_{Q_T}y_\lambda\eta_\lambda d\xi dt\le \int_{Q_T}y\eta d\xi dt.
\end{equation}
 To prove this we first recall equation \eqref{e1.6} which yields  
$$
j_\lambda(y_\lambda+W_G)+j^*_\lambda(\eta_\lambda)=(y_\lambda+W_{G})\eta_\lambda,\quad\mbox{\rm a.e. in}\;Q_T,
$$
and so by    \eqref{e3.6} and since $j_\lambda^*\ge j^*$,  we have
$$
 j(y_\lambda+W_G)+j^* (\eta_\lambda)\le (y_\lambda+W_{G})\eta_\lambda\quad\mbox{\rm a.e. in}\;Q_T,
$$
which yields
$$
\int_{Q_T}( j(y_\lambda+W_G)+j^* (\eta_\lambda)d\xi dt\le \int_{Q_T}(y_\lambda+W_{G})\eta_\lambda d\xi dt.
$$
Since the convex functional
$$
(z,\zeta)\to \int_{Q_T}(j(z)+j^*(\zeta)) d\xi dt
$$
is lower semicontinuous on $L^1(Q_T)$ (and consequently weakly  lower semicontinuous on this space)
we  obtain  
that
\begin{equation}
\label{e3.24}
\int_{Q_T}(j(y+W_G)+j^*(\eta))d\xi dt\le\liminf_{\lambda\to 0}\int_{Q_T}y_\lambda\eta_\lambda d\xi
dt+\int_{Q_T}W_G\eta  d\xi dt.
\end{equation}
Furthermore, by \eqref{e3.12} and again  by the weak lower semicontinuity
of convex integrals in $L^1(Q_T)$ it follows that
\begin{equation}
\label{e3.25}
j(y+W_G),\;j^*(\eta)\in L^1(Q_T).
\end{equation}
On the other hand, since $j(u)+j^*(p)\ge up$ for all $u,p\in \R$ (see \eqref{e1.7}), we have
\begin{equation}
\label{e3.26}
(W_G+y)\eta\le j(y+W_G)+j^*(\eta)\quad\mbox{\rm a. e. in}\;Q_T.
\end{equation}
Moreover, by assumption \eqref{e2.3} we see that for every $M>0$ there exists $R=R(M)\ge 0$, such that
$$
j(-y-W_G)\le M j(y+W_G)\quad\mbox{\rm on}\;Q^R
$$
where
$$
Q^R=\{(t,\xi)\in Q_T:\;|y(t,\xi)+W_G(t,\xi)|\ge R\}.
$$
Since $j(y+W_G)\in L^1(Q_T)$ we have, by continuity of $j$,
\begin{equation}
\label{e3.27}
j(-y-W_G)\le h \quad\mbox{\rm a. e. in}\;Q_T,
\end{equation}
where $h\in L^1(Q_T).$ On the other hand, since $j$ is bounded from below  
 we have
\begin{equation}
\label{e3.27'}
j(-y-W_G) \in L^1(Q_T).
\end{equation}
 Noticing that by virtue of the same
inequality \eqref{e1.7} we have, besides \eqref{e3.26}, that
\begin{equation}
\label{e3.28}
-(y+W_G)\eta\le j(-y-W_G)+j^*(\eta)\quad\mbox{\rm a. e. in}\;Q_T,
\end{equation}
by \eqref{e3.26} and \eqref{e3.27} it follows that a. e. in $Q_T$ we have
$$
|(W_G+y)\eta|\le\max\{j(y+W_G)+j^*(\eta),\;j(-y-W_G)+j^*(\eta)\}\in L^1(Q_T)
$$
and therefore $y\eta\in  L^1(Q_T)$ as claimed (recall that $W_G\in L^\infty(Q_T)$).

Now we come back to equation \eqref{e3.4} which by integration yields
\begin{equation}
\label{e3.29}
\frac12\;\left( |y_\lambda(T)|^2_{-1}-|x|^2_{-1}  \right) +\int_{Q_T}y_\lambda\eta_\lambda d\xi dt+
\lambda\int_{Q_T}y_\lambda(y_\lambda+W_G)d\xi dt=0
\end{equation}
and taking into account that  
\begin{equation}
\label{e3.30}
y_\lambda(T)\to y(T)\quad\mbox{\rm weakly in}\;H, 
\end{equation}
we have by \eqref{e3.29} that
\begin{equation}
\label{e3.31}
\liminf_{\lambda\to 0}\int_{Q_T}y_\lambda\eta_\lambda d\xi dt\le -\frac12\;\left( |y(T)|^2_{-1}-|x|^2_{-1} 
\right).
\end{equation}
In order to complete the proof one needs an integration by parts formula in equation \eqref{e3.20} (or
\eqref{e3.21}) obtained multiplying the equation by $y$ and integrating  on $Q_T$. Formally this is possible
because $y\eta\in \ L^1(Q_T)$ and $y(t)\in H^{-1}(\mathcal O)$ for all $t\in [0,T]$. But, in order to prove it
rigorously, one must give a sense to $(y'(t),y(t))$. Lemma \ref{l3.1} below answers this question
positively and by \eqref{e3.31} also proves \eqref{e3.23}.

We first note that since $j,j^*$ are nonnegative and convex and such that $j(0)=0=j^*(0)$, we have for all measurable
 $f:Q_T\to \R$ that for all $\alpha\in [0,1]$,
 $$j(f)\in L^1(Q_T)\Rightarrow j(\alpha f)\in L^1(Q_T)$$
and
 $$j^*(f)\in L^1(Q_T)\Rightarrow j^*(\alpha f)\in L^1(Q_T).$$
Furthermore, as in the proof of \eqref{e3.27} by  \eqref{e2.3} we obtain
 $$j(f)\in L^1(Q_T)\Rightarrow j(-f)\in L^1(Q_T).$$
By \eqref{e2.3} the latter is, however, also true for $j^*$, if $f\in L^1(Q_T)$ and $\alpha$ is small enough.
Indeed by  \eqref{e2.3} there are $M,R>0$ such that
$$
j(-s)\le Mj(s)\quad\mbox{\rm if}\;|s|\ge R,
$$
hence replacing $s$ by $(-s)$
$$
\frac1M\;j(s)\le j(-s)\quad\mbox{\rm if}\;|s|\ge R.
$$
Now an elementary calculation implies  that for all $p\in \R$
$$
j^*(-p)\le R |p|+\frac1M\;j^*(Mp).
$$
So
$$
j^*(-p/M)\le \frac{R}{M}\; |p|+\frac1M\;j^*(p).
$$
Hence for $\alpha:=1/M$ we have
$$
0\le j^*(-\alpha f)\le \frac{R}{M}\; |f|+\frac1M\;j^*(f)\in L^1(Q_T).
$$
Therefore, $y$ and $\eta$ constructed above fulfill all conditions in the following lemma since
$W_G\in   L^\infty(Q_T)$.
\begin{Lemma}
\label{l3.1}
Let $y\in C^w([0,T];H^{-1}(\mathcal O))\cap L^1(Q_T)$ and   $\eta\in L^1(Q_T)$\break $\cap
L^\infty(0,T;H^{1}(\mathcal O))$
 satisfy
\begin{equation}
\label{e3.31'}
y(t)+A\int_0^t\eta(s)ds=x.
\end{equation}
Furthermore, assume that for some $\alpha>0$, $j(\alpha y), j^*(\alpha \eta)\in L^1(Q_T)$.
Then
$y\eta\in  L^1(Q_T)$,
\begin{equation}
\label{e3.32}
\int_{Q_T}y\eta d\xi dt= -\frac12\;\left( |y(T)|^2_{-1}-|x|^2_{-1} 
\right).
\end{equation}
and
$$
Y_\varepsilon\Sigma_\varepsilon\to y\eta\quad\mbox{\rm in}\;L^1(Q_T),
$$
where $Y_\varepsilon,\Sigma_\varepsilon$ are defined in \eqref{e3.32'} below.
\end{Lemma}
{\bf Proof}. We set for $\varepsilon>0$
\begin{equation}
\label{e3.32'}
Y_\varepsilon=(1+\varepsilon A)^{-m}y,\quad\Sigma_\varepsilon=(1+\varepsilon A)^{-m}\eta,
\end{equation}
where $m\in \N$ is such that $m>\max\{2,(d+2)/4\}$. Then
$$
Y_\varepsilon\in  C^w([0,T];H^{1}_0(\mathcal O)\cap H^{2m-1}(\mathcal O))
\subset  C^w([0,T];H^{1}_0(\mathcal O)\cap C(\overline{\mathcal O}))
$$
and
$$
\Sigma_\varepsilon\in L^1(0,T;W^{2,q}(\mathcal O)),\quad 1<q<\frac{d}{d-1}.
$$
Hence $Y_\varepsilon\Sigma_\varepsilon\in   L^1(Q_T)$ and for $\varepsilon\to 0$
\begin{equation}
\label{e3.33}
\left\{\begin{array}{l}
Y_\varepsilon(t)\to y(t)\quad\mbox{\rm strongly in}\;H^{-1}(\mathcal O),\;\forall\;t\in [0,T]\\
Y_\varepsilon\to y\quad\mbox{\rm strongly in}\;L^1(Q_T)\\
\Sigma_\varepsilon\to \eta\quad\mbox{\rm strongly in}\;L^1(Q_T)\\
\ds\int_0^t\Sigma_\varepsilon(s)ds\to \int_0^t\eta(s)ds\quad\mbox{\rm strongly in}\;H^{1}_0(\mathcal
O)\;\forall\;t\in [0,T].
\end{array}\right.
\end{equation}
We note here that the last fact follows because \eqref{e3.31'} implies that
$\int_0^\bullet\eta(s)ds\in C^w([0,T];H^{1}_0(\mathcal
O)).$
We have also by \eqref{e3.31'}
$$
Y_\varepsilon(t)+A\int_0^t\Sigma_\varepsilon(s)ds=(1+\varepsilon A)^{-m} x,\quad \forall\;t\in [0,T],
$$
which implies
$$
\frac{d}{dt}\;Y_\varepsilon(t)+A\Sigma_\varepsilon(t)=0 
$$
and, taking inner product in $H^{-1}(\mathcal O)$ with $Y_\varepsilon(t)$, we obtain
$$
\frac12\;\frac{d}{dt}\;|Y_\varepsilon(t)|^2_{-1}+\int_\mathcal O\Sigma_\varepsilon(t)Y_\varepsilon(t)d\xi=0,
\quad\mbox{\rm a.e.}\;t\in [0,T].
$$
Hence
\begin{equation}
\label{e3.34}
\lim_{\varepsilon\to 0}\int_{Q_T}\Sigma_\varepsilon(t)Y_\varepsilon(t) d\xi dt=-\frac12\;\left(
|y(T)|^2_{-1}-|x|^2_{-1} 
\right)
\end{equation}
and by \eqref{e3.33} we may assume that for $\varepsilon\to 0$
\begin{equation}
\label{e3.35}
Y_\varepsilon\to y,\quad \Sigma_\varepsilon\to \eta\quad\mbox{\rm a. e. in}\;Q_T.
\end{equation}
Moreover by \eqref{e1.7} we have
\begin{equation}
\label{e3.36}
\alpha^2\Sigma_\varepsilon Y_\varepsilon\le j(\alpha Y_\varepsilon)+j^*(\alpha \Sigma_\varepsilon),\quad
-\alpha^2\Sigma_\varepsilon Y_\varepsilon\le j(-\alpha Y_\varepsilon)+j^*(\alpha \Sigma_\varepsilon)\quad\mbox{\rm a. e.
in}\;Q_T.
\end{equation}
Now we claim  that for $\varepsilon\to 0$
\begin{equation}
\label{e3.37}
j(\alpha Y_\varepsilon)\to j(\alpha y),\;j^*(\alpha \Sigma_\varepsilon)\to j^*(\alpha \eta),\;
j(-\alpha Y_\varepsilon)\to j(-\alpha y)
\quad\mbox{\rm  in}\;L^1(Q_T).
\end{equation}
By \eqref{e3.35}  these convergences hold a.e. in $Q_T$. So, in order to prove \eqref{e3.37}
it suffices to show that $\{j(\alpha Y_\varepsilon)\},
\{j^*(\alpha \Sigma_\varepsilon)\}, \{j(-\alpha Y_\varepsilon)\}$ are equi-integrable on $Q_T$ and
so that they are weakly compact
in 
$L^1(Q_T)$. To this end let $y\in L^1(\mathcal O)$ and let 
$Y_\varepsilon:=(1+\varepsilon A)^{-1}y$, i.e.
$y$ is the solution  to
the equation
\begin{equation}
\label{e3.38zz}
\left\{\begin{array}{l}
Y_\varepsilon-\varepsilon\Delta Y_\varepsilon=y,\quad\mbox{\rm  in}\;\mathcal O,\\
\\
Y_\varepsilon=0,\quad\mbox{\rm  on}\;\partial\mathcal O.
\end{array}\right.
\end{equation}
It may be represented as
\begin{equation}
\label{e3.39zz}
Y_\varepsilon(\xi)=\int_\mathcal OG(\xi,\xi_1)y(\xi_1)d \xi_1,\quad\forall\;\xi\in \mathcal O,
\end{equation}
where $G$ is the associated Green function. It   is well known that $\int_O
G(\xi,\xi_1)d\xi_1$ is the solution to  \eqref{e3.38zz} with $y=1$
so that by the maximum principle we have $0< \int_\mathcal OG(\xi,\xi_1)
d\xi_1 \le 1$
for all $\xi\in \mathcal O$.

We may rewrite $Y_\varepsilon$ as 
$$
Y_\varepsilon(\xi)=\int_\mathcal OG(\xi,\xi_2)d \xi_2\int_\mathcal O\tilde G(\xi,\xi_1)y(\xi_1)d \xi_1
,\quad\forall\;\xi\in \mathcal O,
$$
 where
$$
\tilde G(\xi,\xi_1)=\frac{G(\xi,\xi_1)}{\int_\mathcal OG(\xi,\xi_2)d \xi_2}
$$
and so  $\int_\mathcal O\tilde G(\xi, \xi_1)d \xi_1=1$ for all $\xi\in \mathcal O$.

Then, if $j(y)\in L^1(\mathcal O)$ by Jensen's inequality, since $j(0)=0$ we have
$$
\begin{array}{l}
\ds j(Y_\varepsilon(\xi))\le \int_\mathcal OG(\xi,\xi_2)d  \xi_2\int_\mathcal O\tilde G(\xi, \xi_1)
j(y(\xi_1))d\xi_1
\\
\\
\ds =\int_\mathcal O G(\xi,\xi_1)
j(y( \xi_1))d\xi_1,\quad\forall\;\xi\in \mathcal O.
\end{array}
$$
So, we proved that for any $y\in L^1(\mathcal O)$ with $j(y)\in L^1(\mathcal O)$
$$
j((1+\varepsilon A)^{-1}y)\le (1+\varepsilon A)^{-1}j(y).
$$
Iterating  and using the fact that $(1+\varepsilon A)^{-1}$ preserves positivity we get for  all $m\in \N$
\begin{equation}
\label{e3.40zz}
j((1+\varepsilon A)^{-m}y)\le (1+\varepsilon A)^{-m}j(y), \quad\mbox{\rm a.e. in}\;\mathcal O.
\end{equation}
Now let $y$ be as in the assertion of the lemma and $Y_\varepsilon$ as in \eqref{e3.32'}.
Integrating over $Q_T$, since $(1+\varepsilon A)^{-m}$ is a contraction on $L^1(\mathcal O)$,
\eqref{e3.40zz} applied to $\alpha y$ implies
$$
 \int_{Q_T} j(\alpha Y_\varepsilon(\xi,t))d \xi dt\le\int_{Q_T} j(\alpha y( \xi_2,t))d \xi_2 dt.
$$
Taking into account that $j(\alpha y)\in L^1(Q_T)$ we infer that $\{j(\alpha Y_\varepsilon)\}$ is 
equi-integrable on
$Q_T$. The same argument applies to $\{j^*(\alpha \Sigma_\varepsilon)\}, \{j(-\alpha Y_\varepsilon)\}$.
 
 Then \eqref{e3.36} implies that sequence $\{\Sigma_\varepsilon Y_\varepsilon\}$
is equi-integrable on $Q_T$ and consequently by the Dunford-Pettis theorem, weakly compact in $L^1(Q_T)$.
Since $\{\Sigma_\varepsilon Y_\varepsilon\}$ is a.e. convergent to $y\eta$ we infer that for $\varepsilon\to 0$
\begin{equation}
\label{e3.41zz}
\Sigma_\varepsilon Y_\varepsilon\to y\eta\quad\mbox{\rm strongly in}\;L^1(Q_T),
\end{equation}
which combined with \eqref{e3.34} implies \eqref{e3.32} as desired. $\Box$\bigskip

 We now prove  \eqref{e3.22}.  We have
$$
j(z_\lambda)-j(u)\le \eta_\lambda(z_\lambda-u),\quad \forall\;u\in \R\quad\mbox{\rm a.e. in $Q_T$}. 
$$
Integrating over $Q_T$ yields
$$
\int_{Q_T}j(z_\lambda)d\xi dt\le \int_{Q_T}j(u)d\xi dt+\int_{Q_T}\eta_\lambda(z_\lambda-u)d\xi dt
,\quad \forall\;u\in L^\infty(Q_T).
$$
Note that, by the definition of $\Psi_\lambda$ we have 
$$
z_\lambda=-\lambda\eta_\lambda+y_\lambda+W_G.
$$
Therefore, since $z=y+W_G$, by \eqref{e3.23} and Fatou's lemma we can let
 $\lambda\to 0$ to obtain
$$
\int_{Q_T}j(z)d\xi dt- \int_{Q_T}j(u)d\xi dt\le \int_{Q_T}\eta(z -u)d\xi dt
,\quad \forall\;u\in L^\infty(Q_T).
$$
Now by Lusin's theorem for each $\epsilon>0$ there is a compact subset $Q_\epsilon\subset Q_T$ such that
$(d\xi\otimes dt)(Q_T\setminus Q_\epsilon)\le \epsilon$ and $y,\eta$ are continuous on $Q_\epsilon$. Let $(t_0,x_0)$
be a Lebesgue point
for $y,$ $\eta$ and $y\eta$ and let $B_r$ be  the ball of center $(t_0,x_0)$ and radius $r$. We take
$$
u(t,\xi)=\left\{\begin{array}{l}
z(t,\xi),\quad\mbox{\rm if}\; (t,\xi)\in Q_\epsilon\cap B_r^c\\
v,\quad\mbox{\rm if}\;(t,\xi)\in (Q_\epsilon\cap B_r)\cup (Q_T\setminus Q_\epsilon).
\end{array}\right.
$$
Here $v$ is arbitrary in $\R$. Since $u$ is bounded we can substitute into the above inequality to get
$$
\int_{B_r\cap Q_\epsilon}(j(z)-j(v)-\eta(z-v))d\xi dt\le \int_{Q_T\setminus Q_\epsilon}(\eta(z
-v)+j(v)-j(z))d\xi dt.
$$
Letting $\epsilon\to 0$ we obtain that
$$
\int_{B_r}(j(z)-j(v)-\eta(z-v))d\xi dt\le 0,\quad\forall\;v\in \R,\;r>0.
$$
This yields
$$
j(z(t_0,x_0))\le j(v)+\eta(t_0,x_0)(z(t_0,x_0)-v),\quad\forall\;v\in \R.
$$
and therefore $\eta(t_0,x_0)\in \partial j(z(t_0,x_0))=\Psi(z(t_0,x_0))$.
Since almost all points of $Q_T$ are Lebesgue points we get \eqref{e3.22} as claimed.

 \bigskip

\noindent{\bf Proof of Theorem \ref{t2.3}} (Continued).  Let us first summarize what we have proved for the pair $
(y,\eta)\in L^1(Q_T)\times L^1(Q_T)$. We have
$$
\begin{array}{l}
\ds y\in C^w([0,T];H),\quad \int_0^\bullet\eta(s)ds\in C^w([0,T];H^1_0(\mathcal O)),
\\
\\
\ds  \eta(t,\xi)\in \Psi(y(t,\xi))\quad\mbox{\rm for a.e.}\;(t,\xi)\in Q^T,\\
\\
\ds y(t)+A\int_0^{t}\eta(s)ds=x,\quad t\in [0,T],\\
\\
\ds j(\alpha y),\; j^*(\alpha y)\in L^1(Q_T)\quad\mbox{\rm for some}\;\alpha\in (0,1].
\end{array}
$$
We claim that $(y,\eta)$ is the only such  pair. Indeed,  if $(\tilde y,\tilde \eta)$ is another then
$$
\begin{array}{l}
j(\frac\alpha2(y-\tilde y))\le \frac12\;j(\alpha y)+\frac12\;j(-\alpha\tilde  y)
\end{array}
$$
and
$$
\begin{array}{l}
j^*(\frac\alpha2(y-\tilde y))\le \frac12\;j^*(\alpha y)+\frac12\;j^*(-\alpha\tilde  y).
\end{array}
$$
But as we have explained before Lemma \ref{l3.1} the right hand sides are in $L^1(Q_T)$. Hence
$y-\tilde y$, $ \eta-\tilde\eta$ fulfill all conditions of Lemma \ref{l3.1} and adopting the  notation from there we
have for $\varepsilon>0$
$$
\begin{array}{lll}
\ds Y_\varepsilon(t)-\tilde
Y_\varepsilon(t)&=&\ds\Delta\int_0^t(\Sigma_\varepsilon(s)-\widetilde{\Sigma}_\varepsilon(s))ds
\\
\\
&=&\ds\int_0^t\Delta(\Sigma_\varepsilon(s)-\widetilde{\Sigma}_\varepsilon(s))ds,\quad t\in [0,T].
\end{array}
$$
Differentiating and subsequently taking the inner product in $H$ with $Y_\varepsilon(t)-\tilde
Y_\varepsilon(t)$ and integrating again we arrive at
$$
\begin{array}{l}
\ds\frac12\;\left|(1+\varepsilon A)^{-m}(Y_\varepsilon(t)-\tilde
Y_\varepsilon(t))\right|_{-1}^2=\int_0^t\int_\mathcal O
(Y_\varepsilon(s)-\tilde Y_\varepsilon(s))(\Sigma_\varepsilon(s)-\tilde \Sigma_\varepsilon(s))d\xi ds
\\
\\
\ds=\int_0^t\int_\mathcal O
((1+\varepsilon A)^{-m}(y(s)-\tilde y(s))(1+\varepsilon A)^{-m}(\eta(s)-\tilde
\eta(s))d\xi ds,\quad t\in [0,T].
\end{array}
$$
Letting $\varepsilon\to 0$ and applying Lemma \ref{l3.1} we obtain that for $t\in [0,T]$
$$
\frac12\;|y (t)-\tilde
y (t)|^2_{-1}=
\int_0^t\int_\mathcal O
(y(s)-\tilde y(s)) (\eta(s)-\tilde
\eta(s))d\xi ds\le 0
$$
by the monotonicity of $\Psi$.

Now let us consider the $\omega$-dependence of $y$ and $\eta$.
By \eqref{e3.20}, \eqref{e3.22} we know that $y=y(t,\xi,\omega)$ is the solution to equation
\begin{equation}
\label{e3.42}
\left\{\begin{array}{l}
y'(t)-\Delta\Psi(y(t)+W_G(t)(\omega))=0\quad\mbox{\rm a.e.}\;t\in [0,T],\\
y(0)=x
\end{array}\right. 
\end{equation}
and as seen earlier for $\eta=\eta(t,\xi,\omega)$ as in \eqref{e3.15}
\begin{equation}
\label{e3.43}
\begin{array}{l}
y\in C^w([0,T];H)\cap L^1(Q_T),\quad  \eta\in L^1(Q_T)\\\\
\ds\int_0^\bullet\eta(s)ds\in C^w([0,T];H^1_0(\mathcal O)),
\end{array} 
\end{equation}
and
\begin{equation}
\label{e3.43'}
\eta(t,\xi,\omega)\in \Psi(y(t,\xi,\omega))+W_G(t,\xi,\omega)\quad\mbox{\rm a.e.}\;(t,\xi,\omega)
\in Q_T\times \Omega.
\end{equation}
By the above uniqueness of $(y,\eta)$, it follows that for any sequence $\lambda\to\infty$ we have $\P$-a.s.
$$
\begin{array}{l}
\ds y_\lambda(t)\to y(t)\quad\mbox{\rm weakly in}\;H=H^{-1}(\mathcal O),\;\forall\;t\in [0,T],
\\
\\
\ds y_\lambda \to y \quad\mbox{\rm weakly in}\;L^1(Q_T),\\
\\
\ds\int_0^t\eta_\lambda(s)ds\to\int_0^t\eta(s)ds\quad\mbox{\rm weakly
in}\;L^1(\mathcal O),\;\forall\;t\in [0,T]\;
\\
\mbox{\rm and weakly
in}\;H^{1}_0(\mathcal O),\;\mbox{\rm a.e.}\;t\in [0,T],\\
\\
\ds \eta_\lambda \to \eta \quad\mbox{\rm weakly in}\;L^1(Q_T).
\end{array}
$$
Since $y$ and $\eta$ are hence strong $L^1(Q_T)$-limits of a sequence of convex conbinations of 
$y_\lambda$,  $\eta_\lambda$ respectively, and $y_\lambda$ and  $\eta_\lambda$ are predictable
processes, it follows that so are   $y$ and $\eta$.
In particular, this means that $Y(t)= y(t)+W_G(t)$ is an $H$-valued weakly continuous adapted process and that the
following equation is satisfied
\begin{equation}
\label{e3.44}
Y(t)-\Delta\int_0^t\eta(s)ds=x+\int_0^tG(s)dW(s),\quad t\in [0,T].
\end{equation}
Equivalently
\begin{equation}
\label{e3.45}
\left\{\begin{array}{l}
dY(t)-\Delta \Psi(Y(t))dt=G(t)dW(t),\\
\\
Y(0)=x.
\end{array}\right. 
\end{equation}
In order to prove that $Y$ is a solution  of \eqref{e3.45}
 in the sense of Definition \ref{d2.1} with $G(t)$ replacing
 $B(X(t))$ and to prove uniqueness and some energy estimates for solutions to
equation
\eqref{e3.45} we need an It\^o's formula type result. As in the case of Lemma \ref{l3.1} the difficulty is that the
integral
$$
\int_{Q_T}\Psi(Y)Yd\xi dt
$$ 
might be (in general) not well defined taking into account that  $\Psi(Y),Y\in L^1(Q_T)$ only. We , however, have
\begin{Lemma}
\label{l3.2}
Let $Y$ as above. Then the following equality
holds
\begin{equation}
\label{e3.46}
\begin{array}{l}
\ds\frac12\;|Y(t)|^2_{-1}=\frac12\;|x|^2_{-1}-\int_0^t\int_\mathcal O\eta(s)Y(s)d\xi ds\\
\\
\ds+\int_0^t \langle Y(s),G(s)dW(s)\rangle_{-1}+\frac12\;\int_0^t\|G(s)\|^2_{L_{HS}(L^2(\mathcal O),H)}ds,\quad 
\P\mbox{\rm -a.s}.
 \end{array} 
\end{equation}
Furthermore, $Y\in C_W([0,T];H)\cap L^1((0,T)\times \mathcal O \times  \Omega)$,
and $\eta\in L^1((0,T)\times \mathcal O \times  \Omega)$ and all conditions \eqref{e2.5}-\eqref{e2.8}
are satisfied.
\end{Lemma}
{\bf Proof}.   By Lemma \ref{l3.1} we have that
$Y\eta\in L^1(Q_T)$. Next we introduce the sequences (see the proof of Lemma \ref{l3.1})) for $m\in \N$
$$
Y_\varepsilon=(1+\varepsilon A)^{-m}Y,\quad\Sigma_\varepsilon=(1+\varepsilon A)^{-m}\eta.
$$
For large enough $m$ we can  apply  It\^o's
formula to the problem
\begin{equation}
\label{e3.47}
\left\{\begin{array}{l}
dY_\varepsilon(t)+A\Sigma_\varepsilon(t)=(1+\varepsilon A)^{-m}GdW(t)\\\\
Y_\varepsilon(0)=(1+\varepsilon A)^{-m}x=x_\varepsilon.
\end{array}\right. 
\end{equation}
We have
\begin{equation}
\label{e3.48}
\begin{array}{l}
\ds\frac12\;|Y_\varepsilon(t)|^2_{-1}=\frac12\;|x_\varepsilon|^2_{-1}-\int_0^t\int_\mathcal
O\Sigma_\varepsilon(s)Y_\varepsilon(s)d\xi ds\\
\\
\ds+\int_0^t\langle
Y_\varepsilon(s),G_\varepsilon(s)dW(s)\rangle_{-1}+\frac12\;\int_0^t\|G_\varepsilon(s)\|^2_{L_{HS}(L^2(\mathcal
O),H)}ds,\quad t\in [0,T].
 \end{array} 
\end{equation}
where $G_\varepsilon=(1+\varepsilon A)^{-m}G$. Letting $\varepsilon\to 0$ (since $W_G\in L^\infty(Q_T)$)
  we get by \eqref{e3.41zz}
$$
\int_{Q_T}Y_\varepsilon\Sigma_\varepsilon d\xi ds\to \int_{Q_T}Y \eta  d\xi ds,\quad\P \mbox{\rm a.s.}.
$$ 
Furthermore
$$
Y_\varepsilon(t)\to Y(t)\quad\mbox{\rm strongly in}\;H^{-1}(\mathcal
O),\;\forall\;t\in [0,T],
$$ 
which by virtue of \eqref{e3.48} yields  \eqref{e3.46}, if we can show that for $t\in [0,T]$
\begin{equation}
\label{e3.49b}
\P-\lim_{\varepsilon\to 0}\int_0^t\langle Y_\varepsilon(s),G_\varepsilon(s)dW(s)   \rangle=
\int_0^t\langle Y (s),G (s)dW(s)   \rangle.
\end{equation}
We  shall even show that this convergence in probability is locally uniform in $t$.
We have by a standard consequence of the Burkholder-Davis-Gundy inequality for $p=1$ (see e.g.
\cite[Corollary D-0.2]{14}) that for $\bar Y_\varepsilon:=(1+\varepsilon A)^{-2m}Y$ and $\delta_1,\delta_2>0$
\begin{equation}
\label{e3.50b}
\begin{array}{l}
\ds\P\left[\sup_{t\in[0,T]}\left|\int_0^t\langle Y (s),G (s)dW(s)   \rangle-
\int_0^t\langle Y_\varepsilon(s),G_\varepsilon(s)dW(s) \rangle   \right|\ge \delta_1\right]\\
\\
\ds\le\frac{3 \delta_2}{ \delta_1}+
\P\left[\int_0^T  \|G(s)\|^2_{L_{HS}(L^2(\mathcal O),H)}|Y(s)-\overline{Y}_\varepsilon(s)|^2_{-1}ds    \ge
\delta_2\right].
\end{array} 
\end{equation}
Since $Y\in C^w([0,T];H)$, $\P$-a.s. and $(1+\varepsilon A)^{-1}$ is a contraction on $H$ we have
$$
\sup_{s\in[0,T]}|Y(s)-\overline{Y}_\varepsilon(s)|_{-1}\le 2\sup_{s\in[0,T]}|Y(s)|^2_{-1},\quad\P\mbox{\rm -a.s.}.
$$
Hence by \eqref{e2.13} the second term on the right hand side of \eqref{e3.50b} converges to zero as $\varepsilon\to
0$. Taking subsequently $\delta_2\to 0$, \eqref{e3.50b} implies \eqref{e3.49b}. We emphasize that, since  the left
hand size of \eqref{e3.46} is not continuous $\P$-a.s. (though all terms on the right hand side are), the 
$\P$-zero set of $\omega\in \Omega$ for which \eqref{e3.46} does not hold might depend on $t$.

Next we would like to take expectation in \eqref{e3.46}. Note that because $|Y(t)|_{-1}^2$ is not
$\P$-a.s.continuous in $t$ we cannot use stopping times   to argue that \eqref{e3.46} holds with expectation
taken for every summand and the local martingale term dropped. We need a more delicate argument here. To this end
first note that by \eqref{e3.43'} and \eqref{e1.6} we have
\begin{equation}
\label{e3.50'}
\eta(s)Y(s)=j(Y(s))+j^*(\eta(s))\ge 0,
\end{equation}
hence \eqref{e3.46} implies that for every $t\in [0,T]$
\begin{equation}
\label{e3.51b}
|Y(t)|_{-1}^2\le |x|_{-1}^2+N_t+\int_0^t  \|G(s)\|^2_{L_{HS}(L^2(\mathcal O),H)}ds,\quad\P\mbox{\rm -a.s.},
\end{equation}
where
$$
N_t:=\int_0^t\langle Y (s),G (s)dW(s)   \rangle_{-1},\quad t\ge 0,
$$
is a continuous  local martingale  such that
$$
<N>_t=2\int_0^t  |G^*(s)Y(s)|^2_{L^2(\mathcal O)}ds,\quad t\ge 0,
$$
where $G^*(s)$ is the adjoint of $G(s):L^2(\mathcal O)\to H.$ We shall prove that
\begin{equation}
\label{e3.52b}
\E\left[\sup_{t\in[0,T]}|N_t|\right]<+\infty.
\end{equation}
By the   Burkholder-Davis-Gundy inequality for $p=1$ applied to the stopping times
$$
\tau_N:=\inf\{t\ge 0:\;|N_t|\ge N\}\wedge T,\quad N\in \N,
$$
\begin{equation}
\label{e3.53b}
\begin{array}{lll}
\ds\E\left[\sup_{t\in[0,\tau_N]}|N_t|\right]&\le&\ds 3E
\left[\sup_{s\in[0,\tau_N]}|Y(s)|_{-1}\left(4 \int_0^{\tau_N}  \|G(s)\|^2_{L_{HS}(L^2(\mathcal O);H)}ds 
\right)^{1/2}
\right]\\
\\
&\le&\ds 6C\left(\E\left[\sup_{s\in[0,\tau_N]}|Y(s)|^2_{-1}\right]\right)^{1/2}, 
\end{array}
\end{equation}
where
$$
C:=\left(\E\left[ \int_0^{T}  \|G(s)\|^2_{L_{HS}(L^2(\mathcal O);H)}ds  \right]   \right)^{1/2} <\infty.
$$
Since $Y\in C^w([0,T];H)$, we know that $s\mapsto |Y(s)|^2_{-1}$ is lower semicontinuous. Therefore
by \eqref{e3.51b}
$$
\begin{array}{lll}
\ds\sup_{s\in[0,\tau_N]}|Y(s)|^2_{-1}&=&\ds\sup_{s\in[0,\tau_N]\cap \Q}|Y(s)|^2_{-1}\le |x|^2_{-1}+
\sup_{s\in[0,\tau_N]}|N_s|\\
\\
&&\ds+\int_0^{T}  \|G(s)\|^2_{L_{HS}(L^2(\mathcal O);H)}ds,\quad\P\mbox{\rm -a.s.}.
\end{array}
$$
So \eqref{e3.53b} implies that for all $N\in \N$
$$
\left(\E\left[\sup_{t\in[0,\tau_N]}|N_t|\right]\right)^2\le 36 C^2\left[ |x|^2_{-1}+
\E\left[\sup_{s\in[0,\tau_N]}|N_s|\right]+C^2  \right],
$$
which entails that
$$
\sup_{N\in \N}\E\left[\sup_{t\in[0,\tau_N]}|N_t|\right]<\infty.
$$
By monotone convergence this implies \eqref{e3.52b}, since $N_t$ has continuous sample paths
and $\tau_N\uparrow T$ as $N\to \infty$. Now \eqref{e3.51b} implies that also 
\begin{equation}
\label{e3.54b}
\E\left[\sup_{t\in[0,T]}|Y(t)|^2_{-1}\right]<\infty.
\end{equation}
By \eqref{e3.52b}, \eqref{e3.50'} and \eqref{e3.46} it follows that
\begin{equation}
\label{e3.55b}
\eta Y\in L^1((0,T)\times \mathcal O\times \Omega).
\end{equation}
Hence by \eqref{e3.50'}
$$
j(Y), j^*(\eta)\in L^1((0,T)\times \mathcal O\times \Omega)
$$
and therefore
$$
Y,\eta\in L^1((0,T)\times \mathcal O\times \Omega).
$$
Taking expectation in \eqref{e3.46} we see that $t\mapsto \E[|Y(t)|^2_{-1}]$ is continuous. Since
$Y\in C^w([0,T];H), \P$-a.s., \eqref{e3.54b} then also implies $Y\in C_W([0,T];H)$. This in turn together with
\eqref{e3.44} implies that also \eqref{e2.6} holds. $\Box$

Now we come back to the proof of Theorem \ref{t2.3} noticing that Lemma \ref{l3.2} also
implies  the uniqueness of the solution $Y$ and estimate \eqref{e2.14}.
This concludes the proof of Theorem \ref{t2.3}. $\Box$

\section{Proof of Theorem \ref{t2.2}}
Consider the space
\begin{equation}
\label{e4.1}
\begin{array}{l}
\ds \mathcal K=\Big\{X\in C_W([0,T];H)\cap L^1((0,T)\times\mathcal  O\times \Omega):\;X\;\mbox{\rm predictable},
\\
\\
\hspace{10mm}\ds\sup_{t\in[0,T]}\E [e^{-2\alpha t}|X(t)|^2_{-1}]\le M_1^2,\;\E\int_{Q_T}j(X(s) d\xi ds\le
M_2\Big\},
\end{array}
\end{equation}
where $\alpha >0,M_1>0$ and $M_2>0$ will be specified later.

The space $\mathcal K$ is endowed with the norm
$$
\|X\|_\alpha=\left(\sup_{t\in[0,T]}\E [e^{-2\alpha t}|X(t)|^2_{-1}]\right)^{1/2}.
$$
Note that $\mathcal K$ is closed in the norm $\|\cdot\|_\alpha$. Indeed, if $X_n\to X$ in $\|\cdot\|_\alpha$
then since
$$
\E\int_{Q_T}j(X_n(s))d\xi ds\le M_2,\quad \forall\;n\in \N,
$$
\eqref{e3.14'} implies  that
$$
X_n\to X,\quad \mbox{\rm   in} \;L^1((0,T)\times\mathcal  O\times \Omega)
$$
and by   Fatou's  Lemma we get
$$
\E\int_{Q_T}j(X(s)) d\xi ds\le M_2.
$$
as  claimed.
Now consider the mapping $\Gamma:\mathcal K\to \mathcal K$ defined by
\begin{equation}
\label{e4.2}
Y=\Gamma(X),
\end{equation}
where $Y\in C_W([0,T];H)\cap L^1((0,T)\times\mathcal  O\times \Omega)$ is the solution in the sense of Definition
\ref{d2.1} of the problem
\begin{equation}
\label{e4.3}
\left\{\begin{array}{l}
dY(t)-\Delta\Psi(Y(t))dt=B(X(t))dW(t)\quad\mbox{\rm in}\; Q_T,\\
\Psi(Y(t))=0\quad\mbox{\rm on}\; \Sigma_T,\\
Y(0)=x\quad\mbox{\rm in}\; \mathcal O.
\end{array}\right.
\end{equation}
We shall prove that for $\alpha,M_1,M_2$  suitably chosen,  $\Gamma$ maps $\mathcal K$ into itself and it is a
contraction in the norm $\|\cdot\|_\alpha$.

By \eqref{e4.3}, \eqref{e3.46}  and \eqref{e1.6} we have
$$
\begin{array}{l}
\ds\frac12\;|Y(t)|^2_{-1} +\int_0^t\int_\mathcal O(j(Y(s))+j^*(\eta(s)))d\xi
ds\\
\\
\ds=\int_0^t\langle Y(s),B(X(s))dW(s)\rangle_{-1}\\
\\
\ds+\frac12\;\int_0^t  \|B(X(s))\|^2_{L_{HS}(L^2(\mathcal  O),H)}ds+\frac12\;|x|^2_{-1},\quad t\in [0,T].
 \end{array} 
$$
By Hypothesis $(H_2)$ we have  
$$
\begin{array}{l}
\ds \frac12\;\sup_{t\in[0,T]}\E[e^{-2\alpha t}|Y(t)|^2_{-1}]+2e^{-2\alpha
t}\E\int_0^t\int_\mathcal O (j(Y(s))+j^*(\eta(s)))d\xi ds
\\
\\
\ds \le \frac12\;|x|^2_{-1}  +\frac{L^2}{2}\;\sup_{t\in[0,T]}\left[e^{-2 \alpha t}\int_0^t\E|X(s)|^2_{-1}ds\right]\\
\\
\ds\le \frac12\;|x|^2_{-1} +\frac{L^2}{2}\;\sup_{t\in[0,T]} \int_0^te^{-2\alpha(t-s)}\E e^{-2\alpha
s}|X(s)|^2_{-1}ds\le \frac12\;|x|^2_{-1} + \frac{L^2M_1^2}{4\alpha}.
\end{array}
$$
Hence
$$
\sup_{t\in[0,T]}\E[e^{-2\alpha t}|Y(t)|^2_{-1}]\le \frac{L^2M^2_1}{2\alpha}+|x|_{-1}^2
$$
and
$$\E\int_{Q_T}(j(Y(s))+j^*(\eta(s))))d\xi
\le\left(\frac{L^2M_1^2}{2\alpha}+|x|_{-1}^2\right)e^{2\alpha T}.
$$
Hence for $\alpha>L^2$, $M_1^2>2|x|^2_{-1}$ and  $M_2\ge M_1^2e^{2\alpha T}$   the operator  $\Gamma$ maps
$\mathcal K$ into itself. By a similar computation involving Hypothesis $(H_2)$ we see that
for $M_1, M_2$ and $\alpha$ suitably chosen
\begin{equation}
\label{e4.4}
\|Y_1-Y_2\|_\alpha\le \frac{C}{\sqrt\alpha}\;\|X_1-X_2\|_\alpha
\end{equation}
where $Y_i=\Gamma X_i,\;i=1,2.$ Hence for a suitable $\alpha$, $\Gamma$ is a contraction and so equation
$X=\Gamma(X)$
has a unique solution in $\Gamma$. This completes the proof. $\Box$

\end{document}